\batchmode
\documentclass [11pt,leqno]{article}
\usepackage{amssymb}
\newtheorem{theorem}{Theorem}

\newtheorem{remark}[theorem]{Remark}
\newtheorem{example}[theorem]{Example}

\title{ Singular 2-webs and Polar Curves}
\author{Fernando Etayo\footnote{Departamento de Matem\'{a}ticas, Estad\'{\i}stica y
Computaci\'{o}n.
 Facultad de Ciencias.  Universidad de Cantabria.
 Avda. de los Castros, s/n, 39071 Santander, SPAIN.
 e-mail: etayof@unican.es}}
\date{}
\begin{document}
\maketitle

\begin{abstract}

A 2-web in the plane is given by two everywhere transverse 1-foliations.  In this paper we introduce the study of singular 2-webs, given by any two foliations, which may be tangent in some points.  We show that such two foliations are tangent along a curve, which will be called the polar curve of the 2-web, and we study the relationship between the contact order of leaves of both foliations and the singularities of the polar curve.

\end{abstract}



\section{Introduction}

A 2-web in the real plane is given by two everywhere transverse 1-dimensional foliations, i.e., two families of curves such
that in any point the curves passing through it have different tangent lines. These 2-webs are always locally
diffeomorphic to that of vertical and horizontal lines, and then they have no local invariants. The case of
3-webs is completely different because of the curvature of the Blaschke-Chern connection, which measures how far
the web is to be hexagonal \cite{B, GS}.

A 2-web can be defined by a (1,1) tensor field $F$ of maximum rank such that $F^{2}=I$, where $I$ denotes the
identity tensor field. The eigenspaces associated to the eigenvalues $\pm 1$ define two distributions, which are
involutive in the case of the real plane. $F$ is said to be a paracomplex structure.

One could think that 2-webs have no interest. Nevertheless, we want to go into an unknown landscape, which has
not been explored yet, as far as the author knows. Let us consider two different 1-dimensional foliations in the
plane. What can you say about the set of points in which both foliations are tangent? If this set is non-empty,
we shall say that both foliations define a \emph{singular 2-web} and the set will be called  the \emph{polar
curve} of the singular 2-web. We use this terminology following a similar idea \cite{C} developed in the context
of a holomorphic foliation in $\mathbb{C}^{2}$, where the polar curve is defined from the foliation and a
direction in the plane, i.e., it is defined from 2-web defined by the holomorphic foliation and the foliation of lines parallel to the direction.

Polar curves should not be confused with polar foliations. A polar foliation \cite{A} is  a singular foliation in a complete Riemannian manifold such that for each regular point $p$, there is an immersed submanifold $\Sigma _{p}$, called section, that passes through $p$ and that meets all the leaves and always perpendicularly.

In the previous paper \cite{E2} we have shown some explicit examples. In the present one we state some results
about the polar curve. The following example is introduced as a motivating case.

\begin{example}

Let us consider the punctured plane $\mathbb{R}^{2}-\{(0,0)\}$ and the foliations given by

\begin{itemize}
\item ${\cal F}_{1}$ is the set of circles with center $(0,0)$. \item ${\cal F}_{2}$ is the set of vertical
lines.
\end{itemize}

Obviously, this is a singular 2-web having the horizontal axis as polar curve. We can prove it by using
different techniques: analyzing the associated distributions, the paracomplex structure and the Pfaffian forms. As
is well known 1-dimensional distributions are always integrable, and then we can work interchangeably with
distributions and foliations. We shall show  carefully these three techniques, because we should choose the best
one in order to ahead more complex situations.

\begin{itemize}

\item Associated distributions: The tangent vector of the foliation ${\cal F}_{1}$ at the point $(x,y)$ is
$X=-y\frac{\partial }{\partial x}+x\frac{\partial }{\partial y}$. The tangent vector to ${\cal F}_{2}$ is
$Y=\frac{\partial }{\partial y}$. These vectors are colinear in the polar curve $\{y=0\} -\{(0,0)\}$.

\item The paracomplex structure \cite {EST} is the $(1,1)$-tensor field $F$ such that their eigen\-spaces are
the distributions tangent to the foliations. A straightforward calculation shows that

$$ F= \left( \begin{array}{cc} 1 & 0 \\ -\frac{2x}{y} & -1 \end{array} \right)
$$

One can esaily check that $F^{2}= id; F(X)=X; F(Y)=-Y$. The paracomplex structure $F$ is well defined in all the
punctured plane unless the polar curve $\{y=0\} -\{(0,0)\}$. Then the polar curve appears as the set of points
where the structure cannot be defined.

\item The Pfaffian forms defining ${\cal F}_{1}$ and ${\cal F}_{2}$ are $\omega =x\, dx+y\, dy$ and $\eta = dx$.
The points where these 1-forms are dependent are those where $0=\omega \wedge \eta = (x\, dx+y\, dy)\wedge dx =
-ydx\wedge dy$, which are those of the polar curve.

\end{itemize}
The best option is to work with Pfaffian forms, because one has only to check a product of 1-forms.
\end{example}

\section{Algebraic foliations and paracomplex structure}

As is well known, a foliation is said to be algebraic if it is given by a 1-form $\omega =\omega
_{1}(x,y)dx+\omega _{2}(x,y)dy$, where $\omega _{i}(x,y)$ are polynomials. This doesn't mean that the algebraic
curves of the foliation have to be algebraic curves. For example, the 1-form $\omega =y dx- dy$ is algebraic and
the curves of this foliation are the exponential $C_{k}=\{ y=ke^{x}\}, k\in \mathbb{R}$.

In the same way, we can say that a (1,1)-tensor field $F=\frac{\partial }{\partial x^{i}}\otimes F^{i}_{j}
dx^{j}$ is algebraic (resp. rational) if the functions $F^{i}_{j}$ are polynomial (resp. rational).

\bigskip

Then we have,

\begin{theorem} Let $\omega =\omega _{1}(x,y)dx+\omega _{2}(x,y)dy$ and $\eta =\eta  _{1}(x,y)dx+\eta  _{2}(x,y)dy$
be two 1-forms in an open subset of the real plane. Then,

(1) The polar curve of the singular 2-web defined by $\omega$ and $\eta$ is  the curve $\{ \omega _{1}\eta
_{2}-\omega _{2}\eta  _{1}=0\}$. Besides, if the coefficient functions $\omega _{1},\omega _{2},\eta  _{1},\eta
_{2}$ are polynomial of degrees $p,q,r,s$ then the polar curve is an algebraic curve of degree less or equal to
${\rm max} (p+s, q+r)$.

(2) The paracomplex structure $F$ defined by $\omega$ and $\eta$ is

$$F=\frac{1}{-\omega _{2}\eta  _{1}+\omega _{1}\eta  _{2}}\left( \begin{array}{cc}
-\eta  _{1}\omega _{2}-\omega _{1}\eta  _{2} & -2\omega _{2}\eta  _{2}\\
2\omega _{1}\eta  _{1} & \eta  _{2}\omega _{1} + \omega _{2} \eta  _{1}
\end{array}
\right)
$$

(3) If $\omega$ and $\eta$ are algebraic, then $F$ is rational.

\end{theorem}

\emph{Proof}.

 (1) Observe that $\omega \wedge \eta = ( \omega _{1}\eta _{2}-\omega _{2}\eta  _{1}) dx\wedge dy$,
thus showing that both 1-forms are dependent on the curve $ C= \{ \omega _{1}\eta _{2}-\omega _{2}\eta
_{1}=0\}$.

(2) In the points $p\in \mathbb{R}^{2}-\{ C\}$, where $c$ denotes the polar curve, one has $T_{p}\mathbb{R}^{2}
=(ker\omega )_{p}\oplus (ker\eta )_{p}$. A basis of $ker\omega $ is $-\omega _{2}\frac{\partial }{\partial
x^{1}}+\omega _{1}\frac{\partial }{\partial x^{2}}=(-\omega _{2},\omega _{1})$. A basis of $ker\eta $ is $(-\eta
_{2},\eta _{1})$. Then, a vector field $v=v _{1}\frac{\partial }{\partial x^{1}}+v _{2}\frac{\partial }{\partial
x^{2}}$ can be decomposed as

$$v=\left( \begin{array}{c} v_{1}\\ v_{2} \end{array} \right) =
\alpha \left( \begin{array}{c} -\omega _{2}\\ \omega _{1} \end{array} \right) + \beta \left( \begin{array}{c}
-\eta _{2}\\ \eta _{1}
\end{array} \right)
$$

\noindent which produces, by using the Cramer rule,

$$\alpha = \frac{\mid \begin{array}{cc}v_{1} & -\eta _{2}\\ v_{2} & \eta _{1}\end{array}\mid }
{\mid \begin{array}{cc}-\omega _{2} & -\eta _{2}\\ \omega _{1} & \eta _{1}\end{array}\mid} \hspace{5mm};
\hspace{5mm} \beta = \frac{\mid \begin{array}{cc}-\omega _{2}& v_{1}  \\\omega _{1} & v_{2} \end{array}\mid }
{\mid
\begin{array}{cc}-\omega _{2} & -\eta _{2}\\ \omega _{1} & \eta _{1}\end{array}\mid}
$$

Then, the paracomplex structure $F$ defined as $F|_{ker \omega} =I$, $F|_{ker \eta} =-I$ is that given by

$$F(v)=
\alpha \left( \begin{array}{c} -\omega _{2}\\ \omega _{1} \end{array} \right) - \beta \left( \begin{array}{c}
-\eta _{2}\\ \eta _{1}
\end{array} \right) =
\frac{1}{-\omega _{2}\eta  _{1}+\omega _{1}\eta  _{2}}\left( \begin{array}{cc}
-\eta  _{1}\omega _{2}-\omega _{1}\eta  _{2} & -2\omega _{2}\eta  _{2}\\
2\omega _{1}\eta  _{1} & \eta  _{2}\omega _{1} + \omega _{2} \eta  _{1}
\end{array}
\right) \left( \begin{array}{c} v_{1}\\ v_{2} \end{array} \right)
$$

 (3) It is a direct consequence of (2).

\begin{remark} Observe that the polar curve corresponds to the locus where the paracomplex structure $F$ cannot
be defined. Observe that formula (2) does not depend on the algebricity of foliations.
\end{remark}

\begin{remark} If one changes the 1-forms $\omega $ and $\eta$ by proportional 1-forms $f\omega$ and $g\eta$, f,g beings functions, then the 2-web is the same, and the paracomplex structure $F$ remains invariable in the above theorem. This is important: the 1-forms are not uniquely defined, but the paracomplex structure is uniquely defined, up to sign.

\end{remark}

\section{Singular points of the polar curve}

Let us assume that $\omega $ and $\eta$ are two algebraic foliations with polar curve $C$. As $C$ is an
algebraic curve, $C$ may have singular points. The following examples suggest that singular points of $C$
correspond to the points where the curves of each foliations have contact of order greater than two. Remember
the classical definitions:
\begin{itemize}
\item Two curves $\alpha $ and $\beta$ are said to have \emph{contact of order} $k$ at a point $p$ if their derivatives
of order $0,1,\ldots ,k$ coincide at the point and the derivatives of order $k+1$ are different. We shall denote
it as $ord_{\alpha \beta}(p)=k$.

\item The multiplicity $mult_{C}(p)$ of a curve C at p  is the order of the first non-vanishing term in the
Taylor expansion of f at p, where $C=\{ f(x,y)=0\}$. The point is said to be a \emph{regular} point if $mult_{C}(p)=1$, and \emph{singular} if $mult_{C}(p)\geq 2$.

\item A 1-form $\omega =\omega _{1}(x,y)dx+\omega _{2}(x,y)dy$ is said to define an \emph{exact differential equation} if $\frac{\partial \omega _{1}(x,y) }{\partial y}= \frac{\partial \omega _{2}(x,y) }{\partial x}$. In this case, a curve $f(x,y)=0$ is an integral curve iff $\frac{\partial f}{\partial x}=\omega _{1}(x,y) $ and  $\frac{\partial f}{\partial y}=\omega _{2}(x,y)$. As it is well known, multiplying by an integrating factor $\mu (x,y)$ any 1-form can be transformed into an exact differential equation, although in many cases obtaining the integrating factor is not easy.
\end{itemize}

\begin{example}

 Let us consider the foliations ${\cal F}_{\omega }$ and ${\cal F}_{\eta }$ given by the 1-forms and their dual vector fields:

$\omega =y^{2}dx-dy, \; X_{\omega }=\frac{\partial }{\partial x}+ y^{2}\frac{\partial }{\partial y};$

$\eta =-x^{3}dx+dy, \; X_{\eta }=-\frac{\partial }{\partial x}-x^{3}\frac{\partial }{\partial y}$.

The first one is the set of hyperbolas $\{ y=\frac{1}{k-x}, k\in \mathbb{R}\}$ and the horizontal axis, and the
second one that of quartics $\{ y=k + \frac{x^{4}}{4}, k\in \mathbb{R}\}$. The polar curve is the cubic $C=\{
y^{2}=x^{3}\}$, which has a singular point of multiplicity two (a cusp point) in the origin. The curves of
${\cal F}_{\omega }$ and ${\cal F}_{\eta }$ through the origin are $\{ y=0\}$ and $\{ y=\frac{x^{4}}{4}\}$. The
last one has an inflection point in the origin, thus proving that both curves have contact of order three.

\end{example}

\begin{example} Let us consider the foliation ${\cal F}_{\alpha }$ given by the 1-form $\alpha =dy$, and the
foliation ${\cal F}_{\eta }$ given by the same 1-form $\eta$ of the above example. Then de polar curve is $C=\{
x^{3}=0\}$, which is the vertical axis, being all of its points singular of multiplicity three, equal to the
contact order of tangent curves of each foliation.

\end{example}

We can state:

\begin{theorem} Let ${\cal F}_{\omega }$ and ${\cal F}_{\eta }$ be two algebraic foliations given by two exact differential equations $\omega$ and $\eta$  and let $C$ be their
polar curve. If the curves $\gamma _{\omega }$ and $\gamma _{\eta }$ of
${\cal F}_{\omega }$ and ${\cal F}_{\eta }$ through $p$ have contact of order $\geq 2$ then $p\in c$ is a singular point of the polar curve. In this case, $ord_{\gamma _{\omega } \gamma _{\eta }}(p)\geq mult_{C}(p)$.

\end{theorem}

\emph{Proof}

First we introduce some notation. Let $\omega =\omega _{1}(x,y)dx+\omega _{2}(x,y)dy$ and $\eta =\eta  _{1}(x,y)dx+\eta  _{2}(x,y)dy$ the 1-forms corresponding to ${\cal F}_{\omega }$ and ${\cal F}_{\eta }$ and let

$$C=\{ \omega _{1}\eta
_{2}-\omega _{2}\eta  _{1}=0\} = \{ f(x,y)=0\}$$
\noindent be their polar curve.  Then a point $p\in C$ is singular iff

$$\frac{\partial f}{\partial x}(p) = 0 = \frac{\partial f}{\partial y}(p)
$$
Thus, $p$ is a singular point iff the following equations hold at $p$:

 \begin{equation} \left\{ \begin{array}{cc}
\frac{\partial \omega _{1}}{\partial x} \eta _{2} + \omega _{1} \frac{\partial \eta _{2}}{\partial x}
- \frac{\partial \omega _{2}}{\partial x} \eta _{1} - \omega _{2} \frac{\partial \eta _{1}}{\partial x} & =0\\ & \\
\frac{\partial \omega _{1}}{\partial y} \eta _{2} + \omega _{1} \frac{\partial \eta _{2}}{\partial y}
- \frac{\partial \omega _{2}}{\partial y} \eta _{1} - \omega _{2} \frac{\partial \eta _{1}}{\partial y} & =0
\end{array}
\right\}
\end{equation}

\bigskip

Let $\gamma _{\omega }$ and $\gamma _{\eta }$ the integral curves of ${\cal F}_{\omega }$ and ${\cal F}_{\eta }$ through the point $p$. As the 1-forms $\omega$ and $\eta$ are exact, we have:

 \begin{equation} \left\{ \begin{array}{cc}
\frac{\partial \gamma _{\omega }}{\partial x} & = \omega _{1} \\ & \\
\frac{\partial \gamma _{\omega }}{\partial y} & = \omega _{2}
\end{array}
\right\} \hspace{1cm}
{\rm and}\hspace{1cm}
\left\{ \begin{array}{cc}
\frac{\partial \gamma _{\eta }}{\partial x} & = \eta _{1} \\ & \\
\frac{\partial \gamma _{\eta }}{\partial y} & = \eta _{2}
\end{array}
\right\}
\end{equation}

\bigskip

Then, $ord_{\gamma _{\omega } \gamma _{\eta }}\geq 2$ iff the following equations hold at $p$:

 \begin{equation} \left\{ \begin{array}{ccc}
\frac{\partial \gamma _{\omega }}{\partial x}  = \frac{\partial \gamma _{\eta }}{\partial x} ;  &
\frac{\partial \gamma _{\omega }}{\partial y}  = \frac{\partial \gamma _{\eta }}{\partial y}& \\ & & \\
\frac{\partial ^{2}\gamma _{\omega }}{\partial x^{2}}  = \frac{\partial ^{2}\gamma _{\eta }}{\partial x^{2}} ;  &
\frac{\partial ^{2}\gamma _{\omega }}{\partial x\partial y}  = \frac{\partial ^{2}\gamma _{\eta }}{\partial x\partial y}; &
\frac{\partial ^{2}\gamma _{\omega }}{\partial x^{y}}  = \frac{\partial ^{2}\gamma _{\eta }}{\partial x^{y}}
\end{array}
\right\}
\end{equation}
\bigskip

Then, taking into account equation (2) one easily check that $(3)\Rightarrow (1)$, thus proving the first statement of the theorem. The second one follows from a similar reasoning.

\bigskip

\bigskip

The condition in the above theorem of being $\omega$ and $\eta$ exact differential equations is necessary, as the following example in $\mathbb{R}^{2}-\{ (\pm 1,0)\}$ shows:

\begin{example}

Let  ${\cal F}$ be the vertical axis and the set of circles $C_{a}$, with center in the horizontal axis
passing through the points $(a,0)$ and $(1/a,0)$, where $a\neq \pm 1$.

Points $(a,0), (1/a,0), (1,0),(-1,0)$ define a harmonic quadruple of points, i.e., their cross ratio is -1. The center of $C_{a}$ is the point $\frac{a^{2}+1}{2a}$ and the radius is $\frac{a^{2}-1}{2a}$ thus obtaining the equation

$$ C_{a} =\left\{ \left( x-\frac{a^{2}+1}{2a} \right)^{2}\, +y^{2} \, = \, \left( \frac{a^{2}-1}{2a} \right) ^{2} \right\}= \left\{ x^{2}- \frac{a^{2}+1}{a}\, x +y^{2}+1=0\right\}$$

Let us denote by $f(x,y)=0$  is the equation of the curve $C_{a}$. The tangent line in $p=(x,y)\in C_{a}$ is $f_{x}(x-p_{1})+f_{y}(y-p_{2})$, whose direction is generated by the vector $(-f_{y}(p),f_{x}(p))$, and the corresponding dual form will be $f_{x}(p)dx +f_{y}(p)dy$. In our case,

$$\omega = \left( 2x - \frac{a^{2}+1}{a}\right) \, dx + \, 2y\, dy
$$

Then, for the equation of $C_{a}$, we can deduce that

$$\frac{a^{2}+1}{a} =\frac{x^{2}+y^{2}+1}{x}
$$

\noindent thus showing

$$ \omega =  \frac{x^{2}-y^{2}-1}{x}\, dx \, + \, 2y\, dy \: ; \: {\rm if}\, x\neq 0
$$

We can replace the Pffaf form by another one obtained multiplying by the function $\mu (x,y)=x$, and then we would have coefficients of degree two. Then, we can take

$$ \omega =  (x^{2}-y^{2}-1)\, dx \, + \, 2xy\, dy
$$

\noindent which is not an exact differential equation.

\bigskip

Let us consider the singular 2-web given by the foliations:

\begin{itemize}
\item ${\cal F}_{\omega }$ is the vertical axis and the set of circles $C_{a}$, with center in the horizontal axis
passing through the points $(a,0)$ and $(1/a,0)$. \item ${\cal F}_{\eta }$ is the set vertical lines, whose Pfaffian
form is $\eta =dx$.
\end{itemize}

\noindent Then the polar curve is the reducible curve $\{ xy=0\}$ given by both axis and has a unique singular point. Nevertheless, both foliations have in common the leaf $\{ x=0\}$, and then in all of its points have contact of order infinite, thus showing that there are regular points of the polar curve corresponding to higher order contact of both foliations.

On the other hand, by using Theorem 2 one can obtain the paracomplex structure associated to this 2-web:

$$F=\left( \begin{array}{cc}
1 & 0\\
\frac{-(x^{2}-y^{2}-1)}{xy} & -1
\end{array}
\right)
$$

\noindent thus showing again that the polar curve is $\{ xy=0\}$. But one cannot obtain information about the contact order of leaves of both foliations.

In the present case, as one can easily check, $\mu (x,y)=\frac{1}{x^{2}}$ is an integrating factor of $\omega$, and then we can re-write $\omega$ as

$$ \omega =  \frac{x^{2}-y^{2}-1}{x^{2}}\, dx \, + \, \frac{2y}{x}\, dy
$$

\noindent Then, the equation of the polar curve is $C=\{ -\frac{2y}{x}=0\}$, showing the special property of the vertical axis $\{ x=0\}$. The points where the polar curve cannot be defined when it is obtained from exact differential equations corresponds with those of a common leaf of both foliations. But this example also shows that multiplying the 1-form by an integrating factor can add points where the 1-form is not defined: in the example those of the vertical axis $\{ x=0\}$.
\end{example}

\section{Conclusions}
We write down some global conclusions.
\begin{itemize}
\item The following mathematical objects are equivalent in the real plane: a 2-web, two differential equations, two vector fields, two 1-forms, two distributions. Assume that each of them has no singularities nor zeros (restricting to an open subset of the plane, if necessary). Then, we have two families of curves, and we ask where they are tangent. This is a natural question and, as far as the author knows, there was not yet any answer about it.
\item The points where the foliations are tangent define a curve, called the polar curve of the 2-web. If one defines the 2-web by means of two 1-forms it is very easy to find an equation for the polar curve. If the 1-forms are algebraic, the polar curve is an algebraic curve.
\item The 1-forms associated to a 2-web are not unique. If one takes exact differential equations for them, then we have proved that higher order contact points of the foliations are singular points of the polar curve.
\item There exist integrating factors which allow to obtain an exact differential equation for any 1-form. In general, integrating factors are difficult to be calculated. Besides they can exclude points of the plane.
\item The paracomplex structure is unique up to a sign and the points where it is not defined define the polar curve. One can derive the expression of the paracomplex structure from those of the 1-forms, but there is no a general way for the reverse.
\end{itemize}

\end{document}